\documentclass[11pt]{article}

\usepackage{amssymb}

\makeatletter
\@addtoreset{equation}{section}

\newtheorem{theorem}{Theorem}[section]
\newtheorem{proposition}[theorem]{Proposition}
\newtheorem{corollary}[theorem]{Corollary}
\newtheorem{lemma}[theorem]{Lemma}

\newcommand{\Ecar}[1]{\eta\left(#1\right)+\eta\left(1-#1\right)}
\newcommand{\Eccr}[1]{\eta\left(#1\right)-\eta\left(1+#1\right)}
\newcommand{\tends}[1]{{\displaystyle\mathop{\to}_{#1}}}

\def\endverif{\nopagebreak\newline\mbox{\ }\hfill\rule{2mm}{2mm}}

\def\1{{\bf 1}}

\def\2{{1/2}}

\def\C{{\mathbb C}}
\def\N{{\mathbb N}}
\def\R{{\mathbb R}}
\def\T{{\mathbb T}}
\def\Z{{\mathbb Z}}

\def\A{{\cal A}}
\def\CC{{\cal C}}
\def\P{{\cal P}}
\def\U{{\cal U}}

\def\msigma{\sigma_{-i/2}}
\def\ooplus{\mathop{\oplus}}
\def\notsubset{\hbox{$\subset\hspace{-3.5mm}/\hspace{1.5mm}$}}

\def\CAR{{\rm\underline{CAR}:\ }}
\def\CCR{{\rm\underline{CCR}:\ }}
\def\eps{\varepsilon}

\def\det{{\rm det}}
\def\dim{{\rm dim}\,}
\def\Id{{\rm Id}}
\def\Ip{{\rm Im}}
\def\Ker{{\rm Ker}}
\def\Lin{{\rm Lin}}
\def\Mat{{\rm Mat}}
\def\Rp{{\rm Re}}
\def\Spec{{\rm Spec}\,}
\def\Tr{{\rm Tr}}

\title{\bf ENTROPY OF BOGOLIUBOV AUTOMORPHISMS OF CAR AND CCR ALGEBRAS
WITH RESPECT TO QUASI-FREE STATES}

\author{Sergey V. Neshveyev}

\begin{document}

\maketitle

\begin{abstract}
We compute the dynamical entropy of Bogoliubov automorphisms of CAR and CCR
algebras with respect to arbitrary gauge-invariant quasi-free states. This
completes the research started by St{\o}rmer and Voiculescu, and continued
in works of Narnhofer-Thirring and Park-Shin.
\end{abstract}

\section{Introduction and formulation of main result} \label{1}
One of the most beautiful results in the theory of dynamical entropy is the
formula for the entropy of Bogoliubov automorphisms of the CAR-algebra
with respect to quasi-free states obtained by St{\o}rmer and Voiculescu
\cite{SV} in 1990. They proved it under the assumption that the operator
determining the quasi-free state has pure point spectrum. Since then
several papers devoting to the computation of the entropy of Bogoliubov
automorphisms have appeared. Narnhofer and Thirring \cite{NT2} and Park
and Shin \cite{PS} proved the formula for some operators with continuous
spectrum. The latter paper contains also a similar result for the CCR-algebra.
On the other hand, Bezuglyi and Golodets \cite{BG} proved an analogous
formula for Bogoliubov actions of free abelian groups.

While the cases considered in \cite{NT2} and \cite{PS} required a non-trivial
analysis, the proof of St{\o}rmer and Voiculescu is very elegant. It relies
on an axiomatization of certain entropy functionals on the set of multiplicity
functions. The main axiom there stems from the equality
$h_\omega(\alpha)={1\over n}h_\omega(\alpha^n)$. Thus, their method can not
be directly applied to groups without finite-index subgroups. Instead,
we can "cut and move" multiplicity functions without changing the entropy
(see Lemma~\ref{5.1} below). This observation together with the methods
developed in \cite{BG} allowed to prove (under the same restrictions on
quasi-free states) an analogue of
St{\o}rmer-Voiculescu's formula for Bogoliubov actions of arbitrary
torsion-free abelian groups \cite{GN2}. In this paper we will show that,
in fact, the formula holds without any restrictions on the operator
determining the quasi-free state. We will prove also an analogous result
for the CCR-algebra.

We will consider only the case of single automorphism, since in view of
the methods of \cite{GN2} the case of arbitrary torsion-free abelian group
gives nothing but more complicated notations. So the main result of the
paper is as follows.

\begin{theorem} \label{1.1}
 Let $U$ be a unitary operator on a Hilbert space $H$, $\alpha_U$ the
 corresponding Bogoliubov automorphism of the CAR or the CCR algebra over $H$,
 $A$ a bounded ($A\le1$ for CAR) positive operator commuting with $U$ and
 determining a quasi-free
 state $\omega_A$. Let $U_a=U|_{H_a}$ be the absolutely continuous part
 of $U$,
 $$
 H_a=\int^\oplus_\T H_zd\lambda(z),\ \ U_a=\int^\oplus_\T z\,d\lambda(z),\ \
  A|_{H_a}=\int^\oplus_\T A_zd\lambda(z)
 $$
 a direct integral decomposition, where $\lambda$ is the Lebesgue measure
 on $\T$ ($\lambda(\T)=1$). Then

 \CAR $\displaystyle h_{\omega_A}(\alpha_U)=\int_\T
    \Tr(\Ecar{A_z})d\lambda(z)$,

 \CCR $\displaystyle h_{\omega_A}(\alpha_U)=\int_\T
    \Tr(\Eccr{A_z})d\lambda(z)$.
\end{theorem}

\begin{corollary} \label{1.2}
 The necessary condition for the finiteness of the entropy is that $A_z$ has
 pure point spectrum for almost all $z\in\T$.
\end{corollary}

\begin{corollary} \label{1.3}
 If the spectrum of the unitary operator is singular, then the entropy is
 zero.
\end{corollary}

For CAR, the latter corollary is already known from \cite{SV}.

Finally, for systems considered in \cite{NT2} and \cite{PS}, Theorem \ref{1.1}
may be reformulated as

\begin{corollary} \label{1.4}
 Let $I$ be an open subset of $\R$, $\omega$ a locally absolutely continuous
 function
 on $I$, $\rho$ a bounded ($\rho\le1$ for CAR) positive measurable function
 on $I$. Let $U$ and $A$ be the operators on $L^2(I,dx)$ of multiplication
 by the functions $e^{i\omega}$ and $\rho$, respectively. Then

 \CAR $\displaystyle h_{\omega_A}(\alpha_U)={1\over2\pi}\int_I
    [\Ecar{\rho(x)}]|\omega'(x)|dx$,

 \CCR $\displaystyle h_{\omega_A}(\alpha_U)={1\over2\pi}\int_I
    [\Eccr{\rho(x)}]|\omega'(x)|dx$.
\end{corollary}

\smallskip
The paper is organized as follows. Section \ref{2} contains some
preliminaries on entropy and algebras of canonical commutation and
anti-commutation relations. In Section \ref{3} we prove that the entropies
don't exceed the values of the integrals in Theorem \ref{1.1}. The opposite
inequality
is proved in Sections \ref{4} and \ref{5}. In Section \ref{4} we obtain
a lower bound for the entropy in the case where the unitary operator has
Lebesgue spectrum and the operator determining the quasi-free state is close
to a scalar operator. In Section \ref{5}, first, using the observation
mentioned above we extend the estimate of Section \ref{4} to arbitrary
unitaries, and then prove the required inequality.

There are also two appendices to the paper. The results of \cite{GN1} show
that modular automorphisms can have the K-property (in the sense
of Narnhofer and Thirring \cite{NT1}). This observation combined
with the results of the present paper allow to construct on the
hyperfinite III$_1$-factor a simple example of non-conjugate
K-systems with the same finite entropy. This is done in
Appendix~A. Appendix~B contains an auxiliary result on
decomposable operators.

\bigskip\bigskip
\section{Preliminaries} \label{2}
Recall the definition of dynamical entropy \cite{CNT}. Let $(A,\phi,\alpha)$
be a C$^*$-dynamical system, where $A$ is a C$^*$-algebra, $\phi$ a state on
$A$, $\alpha$ a $\phi$-preserving automorphism of~$A$. By a channel in $A$ we
mean a unital completely positive mapping $\gamma\colon B\to A$ of a
finite-dimensional C$^*$-algebra~$B$. The mutual entropy of channels
$\gamma_i\colon B_i\to A$, $i=1,\ldots,n$, with respect to $\phi$ is given
by
$$
H_\phi(\gamma_1,\ldots,\gamma_n)=\sup\sum_{i_1,\ldots,i_n}
   \eta(\phi_{i_1\ldots i_n}(1))+\sum^n_{k=1}\sum_{i_k}S(\phi\circ\gamma_k,
   \phi^{(k)}_{i_k}\circ\gamma_k),
$$
where $\eta(t)=-t\log t$, $S(\cdot\,,\,\cdot)$ the relative entropy,
$\phi^{(k)}_{i_k}=\sum_{i_1,\ldots,\hat i_k,\ldots,i_n}\phi_{i_1\ldots i_n}$,
and the supremum is taken over all finite decompositions
$\phi=\sum\phi_{i_1\ldots i_n}$ of $\phi$ in the sum of positive linear
functionals. If $A$ is a W$^*$-algebra and $\phi$ is a normal faithful state,
then any positive linear functional $\psi\le\phi$ on $A$ is of the form
$\phi(\cdot\,\msigma^\phi(x))$ for some $x\in A$, $0\le x\le1$. Thus,
$$
H_\phi(\gamma_1,\ldots,\gamma_n)=\sup\sum_{i_1,\ldots,i_n}
   \eta(\phi(x_{i_1\ldots i_n}))+\sum^n_{k=1}\sum_{i_k}S(\phi(\gamma_k(\cdot)),
   \phi(\gamma_k(\cdot)\msigma(x^{(k)}_{i_k}))),
$$
where the supremum is taken over all finite partitions of unit.

The entropy of the automorphism $\alpha$ with respect to a channel $\gamma$
and the state $\phi$ is given by
$$
h_\phi(\gamma;\alpha)=\lim_{n\to\infty}{1\over n}
   H_\phi(\gamma,\alpha\circ\gamma,\ldots,\alpha^{n-1}\circ\gamma).
$$
The entropy $h_\phi(\alpha)$ of the system $(A,\phi,\alpha)$ is the supremum
of $h_\phi(\gamma;\alpha)$ over all channels $\gamma$ in~$A$.

We refer the reader to \cite{CNT}, \cite{OP}, \cite{SV}, \cite{NT1} for
general properties of entropy.

\begin{lemma} \label{2.1}
 Let $(A,\phi,\alpha)$ be a C$^*$-dynamical system, $\{A_n\}^\infty_{n=1}$
 a sequence of $\alpha$-invariant subalgebras of $A$, $\{F_n\}^\infty_{n=1}$
 a sequence of completely positive unital mappings $F_n\colon A\to A_n$ such
 that $||F_n(x)-x||_\phi\to0$ as $n\to\infty$, for any $x\in A$. Then
 $$
 h_\phi(\alpha)\le\liminf_{n\to\infty}h_\phi(\alpha|_{A_n}).
 $$
\end{lemma}

\noindent{\it Proof.} The result follows from the continuity of mutual entropy
in $||\ ||_\phi$-topology: see the proof of Lemma~3.3 in \cite{SV}.
\endverif

Though the possibility of $A_n\notsubset A_{n+1}$ is important for
applications to actions of more general groups (see the proof of Theorem 4.1
in \cite{GN2}), we will use this lemma only when $A_n\subset A_{n+1}$. Then
the existence of $F_n$'s is not necessary, as the following result shows.

\begin{lemma} \label{2.2}
 Let $A$ be a C$^*$-algebra, $\phi$ a state on $A$, $\{A_n\}^\infty_{n=1}$
 an increasing sequence of C$^*$-subalgebras such that $\cup_n\pi_\phi(A_n)$
 is weakly dense in $\pi_\phi(A)$. Then, for any channel
 $\gamma\colon B\to A$ and any $\eps>0$, there exist $n\in\N$ and a channel
 $\tilde\gamma\colon B\to A_n$ such that $||\gamma-\tilde\gamma||_\phi<\eps$.
\end{lemma}

\noindent{\it Proof.} This follows from the identification of completely
positive maps $\Mat_d(\C)\to A$ with positive elements in $\Mat_d(A)$
\cite{CE} and, in fact, is implicitly contained in \cite{CNT}. We include
a proof for the convenience of the reader.

Without loss of generality we may suppose that $B=\Mat_d(\C)$. The channels
$B\to A$ are in one-to-one correspondence with positive elements
$Q\in\Mat_d(A)$ such that $\sum_kQ_{kk}=1$. By Kaplansky's density theorem,
there exists a net $\{\tilde Q_i\}_i\subset\cup_n\Mat_d(\pi_\phi(A_n))$ such
that
$$
0\le\tilde Q_i\le1,\ \ \tilde Q_i\tends i\pi_\phi(Q) \ \ \hbox{strongly}.
$$
We can lift $\tilde Q_i$ to an element $Q_i\in\cup_n\Mat_d(A_n)$,
$0\le Q\le1$. For $\delta>0$, set
$$
Q(i;\delta)_{kl}=\left(\sum_jQ(i)_{jj}+d\delta\right)^{-\2}
   (Q(i)_{kl}+\delta_{kl}\delta)\left(\sum_jQ(i)_{jj}+d\delta\right)^{-\2}.
$$
Let $\gamma_{i,\delta}\colon B\to\cup_nA_n$ be the corresponding channel,
$\gamma_{i,\delta}(e_{kl})=Q(i;\delta)_{kl}$. Then
$$
\displaystyle\lim_{\delta\to0}\lim_i||\gamma-\gamma_{i,\delta}||_\phi=0.
$$
\endverif

\medskip
Now recall some facts concerning CAR and CCR algebras \cite{BR2}.

Let $H$ be a Hilbert space. The CAR-algebra $\A(H)$ over $H$ is a
C$^*$-algebra generated by elements $a(f)$ and $a^*(f)$, $f\in H$, such that
the mapping $f\mapsto a^*(f)$ is linear, $a(f)^*=a^*(f)$ and
$$
a^*(f)a(g)+a(g)a^*(f)=(f,g)1, \ \ a(f)a(g)+a(g)a(f)=0.
$$

Each unitary operator $U$ on $H$ defines a Bogoliubov automorphism
$\alpha_U$ of $\A(H)$, $\alpha_U(a(f))=a(Uf)$. The fixed point algebra
$\A(H)_e=\A(H)^{\alpha_{-1}}$ is called the even part of $\A(H)$.

Each operator $A$ on $H$, $0\le A\le1$, defines a quasi-free state $\omega_A$
on $\A(H)$,
$$
\omega_A(a^*(f_1)\ldots a^*(f_n)a(g_m)\ldots a(g_1))=\delta_{nm}
  \det((Af_i,g_j))_{i,j}.
$$

If $\Ker\,A=\Ker(1-A)=0$, then $\omega_A$ is a KMS-state,
\begin{equation} \label{e2.1}
\sigma^{\omega_A}_t(a(f))=a(B^{it}f),\ \  \hbox{where}\ \
 B={A\over1-A}.
\end{equation}

If $U$ and $A$ commute, then $\omega_A$ is $\alpha_U$-invariant.

If $H=K\oplus L$, then $\A(K)$ and $\A(L)_e$ commute, and we have
$$
\A(H)^{\alpha_{1\oplus -1}}=\A(K)\vee\A(L)_e\cong\A(K)\otimes\A(L)_e.
$$
If $K$ is an invariant subspace for $A$, then
$$
\omega_A|_{\A(K)\otimes\A(L)_e}=\omega_A|_{\A(K)}\otimes\omega_A|_{\A(L)_e}.
$$
In particular, there exists an $\omega_A$-preserving conditional expectation
$$
\left(\Id_{\A(K)}\otimes\omega_A(\cdot)|_{\A(L)_e}\right)\circ
   {1+\alpha_{1\oplus -1}\over2}
$$
onto $\A(K)$. If $K$ is finite-dimensional, then $\A(K)$ is a full matrix
algebra of dimension $2^{2n}$. In particular, for any $f\in H$, $||f||=1$,
the algebra $\A(\C f)$ is isomorphic to $\Mat_2(\C)$, and we define matrix
units for it as
\begin{equation} \label{e2.2}
e_{11}(f)=a(f)a^*(f), \ e_{22}(f)=a^*(f)a(f),\ e_{12}(f)=a(f),\
 e_{21}(f)=a^*(f).
\end{equation}
The restriction of a quasi-free state $\omega_A$ to $\A(\C f)$ is given
by the matrix
\begin{equation} \label{e2.3}
\pmatrix{1-\lambda & 0\cr
         0 & \lambda},\ \ \hbox{where}\ \ \lambda=(Af,f).
\end{equation}

\smallskip
The CCR-algebra $\U(H)$ over $H$ is a C$^*$-algebra generated by unitaries
$W(f)$, $f\in H$, such that
$$
W(f)W(g)=e^{i{\Ip(f,g)\over2}}W(f+g).
$$
A representation $\pi$ of $\U(H)$ is called regular, if the mapping
$\R\ni t\mapsto\pi(W(tf))$ is strongly continuous. For any such a
representation, the generator $\Phi_\pi(f)$ of the group $\{\pi(W(tf))\}_t$
is defined, $\pi(W(tf))=e^{it\Phi_\pi(f)}$. Then annihilation and creation
operators are defined as
$$
a_\pi(f)={\Phi_\pi(f)+i\Phi_\pi(if)\over\sqrt2},\ \
 a^*_\pi(f)={\Phi_\pi(f)-i\Phi_\pi(if)\over\sqrt2}.
$$
These are closed unbounded operators affiliated with $\pi(\U(H))''$,
$a_\pi(f)^*=a^*_\pi(f)$, $a^*_\pi(f)$ depends on $f$ linearly, and for any
$f,g\in H$ we have the commutation relations
$$
a_\pi(g)a^*_\pi(f)-a^*_\pi(f)a_\pi(g)=(f,g)1, \ \
 a_\pi(g)a_\pi(f)-a_\pi(f)a_\pi(g)=0
$$
on a dense subspace. In the sequel we will suppress $\pi$ in the notations
of annihilation and creation operators.

Each unitary operator $U$ on $H$ defines a Bogoliubov automorphism
$\alpha_U$ of $\U(H)$, $\alpha_U(W(f))=W(Uf)$.

Each positive operator $A$ on $H$ defines a quasi-free state $\omega_A$
on $\U(H)$,
$$
\omega_A(W(f))=e^{-{1\over4}||f||^2-{1\over2}(Af,f)}.
$$
The cyclic vector $\xi_{\omega_A}$ in the GNS-representation belongs to
the domain of any operator of the form $a^\#(f_1)\ldots a^\#(f_n)$, where
$a^\#$ means either $a^*$ or $a$, and
$$
(a^*(f)a(g)\xi_{\omega_A},\xi_{\omega_A})=(Af,g).
$$
If $\Ker\,A=0$, then $\omega_A$ is separating (i.~e., $\xi_{\omega_A}$ is
separating for $\pi_{\omega_A}(\U(H))''$), and
$$
\sigma^{\omega_A}_t(W(f))=W(B^{it}f),\ \  \hbox{where}\ \
 B={A\over1+A},
$$
so that
\begin{equation} \label{e2.4}
\Delta^{it}_{\omega_A}a^\#(f_1)\ldots a^\#(f_n)\xi_{\omega_A}
 =a^\#(B^{it}f_1)\ldots a^\#(B^{it}f_n)\xi_{\omega_A}.
\end{equation}

If $H=K\oplus L$, then $\U(H)\cong\U(K)\otimes\U(L)$.
If $K$ is an invariant subspace for $A$, then
$$
\omega_A=\omega_A|_{\U(K)}\otimes\omega_A|_{\U(L)},
$$
so that there exists an $\omega_A$-preserving conditional expectation
$\Id_{\U(K)}\otimes\omega_A|_{\U(L)}$
onto $\U(K)$.

If $K$ is finite-dimensional, then every regular representation
$\pi$ of $\U(K)$ is quasi-equivalent to the Fock representation, in
particular, $\pi(\U(K))''$ is a factor of type I$_\infty$ (if $K\ne0$).
Thus, for any regular state $\omega$ on $\U(K)$ (so that the mapping
$t\mapsto\omega(W(tf))$ is continuous) the von Neumann entropy of the
continuation $\bar\omega$ of the state $\omega$ to $\pi_\omega(\U(K))''$ is
defined. We will denote it by $S(\omega)$ (in fact, the notion of entropy of
state
can be defined for all C$^*$-algebras, and then $S(\omega)=S(\bar\omega)$
\cite{OP}). If $K=\C f$, $||f||=1$, we define a system of matrix units
$\{e_{ij}(f)\}_{i,j\in\Z_+}$ for $\pi(\U(K))''$ as follows:
\begin{quote}
$e_{kk}(f)$ is the spectral projection of $a^*(f)a(f)$ corresponding to
$\{k\}$,
\end{quote}
\begin{equation} \label{e2.5}
e_{k+n,k}(f)=\left({k!\over(k+n)!}\right)^\2a^*(f)^ne_{kk}(f)
 =\left({k!\over(k+n)!}\right)^\2\overline{e_{k+n,k+n}(f)a^*(f)^n}.
\end{equation}
In particular, if $\omega_A$ is a quasi-free state on $\U(H)$, for any
$f\in H$, $||f||=1$, we obtain a system of matrix units
$\{e_{ij}(f)\}_{i,j}$ in $\pi_{\omega_A}(\U(H))''$, and
\begin{equation} \label{e2.6}
\omega_A(e_{ij}(f))=\delta_{ij}{\lambda^i\over(1+\lambda)^{i+1}},
 \ \ \hbox{where} \ \ \lambda=(Af,f).
\end{equation}
(This is equivalent to the fact that if $A$ is of trace class, then the
quasi-free state $\omega_A$ is given in the Fock representation by the
density operator ${\Gamma(B)\over\Tr\,\Gamma(B)}$, where $\Gamma$
is the operator of second quantization.)

\smallskip
In the sequel we will write $\CC(H)$ instead of $\A(H)$ and $\U(H)$ in the
arguments that are identical for CAR and CCR.

The following result is known, but we will give a proof for the reader's
convenience.

\begin{lemma} \label{2.3}
 Let H be finite-dimensional, $\omega_A$ a quasi-free state on $\CC(H)$.
 Then

 (i) \CAR $S(\omega_A)=\Tr(\Ecar A)$,\ \CCR $S(\omega_A)=\Tr(\Eccr A)$;

 (ii) if $H=H_1\oplus H_2$, then
 $S(\omega_A)\le S(\omega_A|_{\CC(H_1)})+S(\omega_A|_{\CC(H_2)})$.
\end{lemma}

\noindent{\it Proof.} Let $P_i$ be the projection onto $H_i$,
$A_i=P_iA|_{H_i}$. Set $M_i=\pi_{\omega_{A_i}}(\U(H_i))''$,
$M=\pi_{\omega_A}(\U(H))''$. Since all regular representations of $\U(H_i)$
are quasi-equivalent, we may consider $M_i$ as a subalgebra of $M$. Since
$M_1$ is a type I factor, we have $M=M_1\otimes(M_1'\cap M)$, whence
$M=M_1\otimes M_2$. Thus, the assertion (ii) for CCR is the usual
subadditivity of von Neumann entropy.

Turning to CAR, let us first note
that if $M$ is a full matrix algebra, $\omega$ a state on $M$ and $\alpha$
an automorphism of $M$, then $S(\omega)\le S(\omega|_{M^\alpha})$, and the
equality holds iff $\omega$ is $\alpha$-invariant. Indeed, let $Q$ (resp.
$\tilde Q$) be the density operator for $\omega$ (resp. $\omega|_{M^\alpha}$).
Since the canonical trace on $M^\alpha$ is given by the restriction of the
canonical trace $\Tr$ on $M$, we have $\Tr\,\tilde Q=1$, hence
$$
S(\omega|_{M^\alpha})-S(\omega)=\Tr\,Q(\log Q-\log\tilde Q)\ge0,
$$
and the equality holds iff $Q=\tilde Q$, i.~e., $Q\in M^\alpha$.

Applying this to CAR, we obtain
$$
S(\omega_A)\le S(\omega_A|_{\A(H_1)\otimes\A(H_2)_e})
 \le S(\omega_A|_{\A(H_1)})+S(\omega_A|_{\A(H_2)_e})
 =S(\omega_A|_{\A(H_1)})+S(\omega_A|_{\A(H_2)}).
$$

We see also that if $H_i$ is an invariant subspace for $A$, then
$$
S(\omega_A)=S(\omega_A|_{\CC(H_1)})+S(\omega_A|_{\CC(H_2)}).
$$
So, in proving (i) it is enough to consider one-dimensional spaces, for
which the result follows immediately from (\ref{e2.3}) and (\ref{e2.6}).
\endverif

\begin{lemma} \label{2.4}
 Let $U$ be a unitary operator on $H$, $\{P_n\}^\infty_{n=1}$ a sequence of
 projections in $B(H)$, $P_nU=UP_n$, $P_n\to1$ strongly, $H_n=P_nH$. Then, for
 the Bogoliubov automorphism $\alpha_U$ and any $\alpha_U$-invariant
 quasi-free state $\omega_A$ on $\CC(H)$, we have
 $$
 h_{\omega_A}(\alpha_U)\le\liminf_{n\to\infty}
    h_{\omega_A}(\alpha_U|_{\CC(H_n)}).
 $$
\end{lemma}

\noindent{\it Proof.} Let $C$ be an operator commuting with $P_n$ for all
$n\in\N$. Let $E_n$ be the $\omega_C$-preserving conditional expectation
of $\CC(H)$ onto $\CC(H_n)$ defined above. Then
$||E_n(x)-x||_{\omega_A}\to0$ for any $x\in\CC(H)$. Indeed, for CAR we have
even the convergence in norm, that follows from $||E_n||=1$ and
$||a(f)||=||f||$. For CCR, the assertion follows from the equalities
$$
E_n(W(f))=e^{-{1\over4}||(1-P_n)f||^2-{1\over2}(C(1-P_n)f,(1-P_n)f)}W(P_nf),
$$
$$
||W(f)-W(g)||^2_{\omega_A}=2-2\Rp\left(e^{i{\Ip(f,g)\over2}}
   \omega_A(W(f-g))\right).
$$
Thus we can apply Lemma~\ref{2.1}.
\endverif

\bigskip\bigskip
\section{Upper bound for the entropy} \label{3}
In this section we will prove that the entropies do not exceed the values
of the integrals in Theorem \ref{1.1}.

There exists a Hilbert space $K$ and a unitary operator $V$ on $K$ such that
$U_a\oplus V$ has countably multiple Lebesgue spectrum. Set
$$
\tilde H=H\oplus K,\ \ \tilde U=U\oplus V,\ \ \tilde A=A\oplus0.
$$
Then, due to the existence of an $\omega_{\tilde A}$-preserving conditional
expectation $\CC(\tilde H)\to\CC(H)$, we have
$h_{\omega_A}(\alpha_U)\le h_{\omega_{\tilde A}}(\alpha_{\tilde U})$. On the
other hand, the passage to $(\tilde H,\tilde U,\tilde A)$ does not change
the value of the integral in Theorem \ref{1.1}. So, without loss of
generality we may
suppose that $U_a$ has countably multiple Lebesgue spectrum. If the value
of the integral
is finite, then $A_z$ has pure point spectrum for almost all $z\in\T$. Then
we can represent $H_a$ as the sum of a countable set of copies of $L^2(\T)$
in such a way that $U$ and $A$ act on the $n$-th copy as multiplications by
functions $z$ and $\lambda_n(z)$, respectively (see Appendix~B).
By Lemma~\ref{2.4}, we may restrict ourselves to the sum of a finite number
of copies of $L^2(\T)$. Thus, we suppose
$$
H_a=\ooplus^{m_0}_{k=1}L^2(\T),\ \ U_a=\ooplus^{m_0}_{k=1}z,\ \
 A|_{H_a}=\ooplus^{m_0}_{k=1}\lambda_k(z),
$$
and we have to prove that

\CAR $\displaystyle h_{\omega_A}(\alpha_U)\le\sum^{m_0}_{k=1}
   \int_\T(\Ecar{\lambda_k(z)})d\lambda(z)$,

\CCR $\displaystyle h_{\omega_A}(\alpha_U)\le\sum^{m_0}_{k=1}
   \int_\T(\Eccr{\lambda_k(z)})d\lambda(z)$.

Let $H_0$ be the $m_0$-dimensional subspace of $H$ spanned by constant
functions in each copy of $L^2(\T)$. Then $H_a=\oplus_{n\in\Z}U^nH_0$.
For $n\in\N$, set $H_n=\oplus^n_{k=0}U^kH_0$. We state that
\begin{equation} \label{e3.1}
h_{\omega_A}(\alpha_U)
 \le\lim_{n\to\infty}{1\over n}S(\omega_A|_{\CC(H_{n-1})}).
\end{equation}
For CAR, this is implicitly contained in the proof of Lemma~5.3 in \cite{SV}.
So we will consider CCR only.

For a finite set $X$, we denote by $\Mat(X)$ the C$^*$-algebra of linear
operators on $l^2(X)$. Let $\{e_{xy}\}_{x,y\in X}$ be the canonical system
of matrix units for $\Mat(X)$. Following Voiculescu (see Lemmas 5.1 and 6.1
in \cite{V}), for
$X\subset H$, we introduce unital completely positive mappings
$$
i_X\colon\Mat(X)\to\U(H),\ \ j_X\colon\U(H)\to\Mat(X),
$$
$$
i_X(e_{xy})={1\over|X|}W(x)W(y)^*, \ \
 j_X(a)=P_X\pi_\tau(a)P_X,
$$
where $\tau$ denotes the unique trace on $\U(H)$ ($\tau(W(f))=0$ for $f\ne0$),
and $P_X$ is the projection onto the subspace
$\Lin\{\pi_\tau(W(x))\xi_\tau\,|\,x\in X\}\subset H_\tau$
identified with $l^2(X)$. Then
$$
(i_X\circ j_X)(W(f))={|X\cap(X-f)|\over|X|}W(f)\ \ \forall f\in H.
$$
Hence, for any subspace $K$ of $H$, there exists a net $\{X_i\}_i$ of finite
subsets of $K$ such that
$||(i_{X_i}\circ j_{X_i})(a)-a||\tends{i}0$ $\forall a\in\U(K)$.

Let $H_s=H\ominus H_a$ be the subspace corresponding to the singular part of
the spectrum of $U$. By Lemma~\ref{2.2}, in computing the entropy we may
consider only the channels in $\cup_m\U(H_s\oplus H_m)$. If $\gamma$ is a
channel in $\U(H_s\oplus H_m)=\U(H_s)\otimes\U(H_m)$, then it can be
approximated in norm by a
channel of the form $(i_X\otimes i_Z)\circ(j_X\otimes j_Z)\circ\gamma$,
where $X\subset H_s$ and $Z\subset H_m$. Hence, it suffices to consider
only the channels $i_X\otimes i_Z$.

So, let
$\gamma=i_X\otimes i_Z\colon\Mat(X)\otimes\Mat(Z)\to
  \U(H_s)\otimes\U(H_m)=\U(H_s\oplus H_m)$. Set $L=\Lin\,X$. Fix $\eps>0$.
By Lemma~5.1 in \cite{SV}, there exist $n_0\in\N$ and a sequence of
projections $\{Q_n\}^\infty_{n=n_0}$ in $B(H_s)$ such that
$\dim Q_n\le\eps n$ and $||(U^k-Q_nU^k)|_L||\le\eps$ for $k=0,\ldots,n-1$.
Define a channel $i^{(n,k)}_X\colon\Mat(X)\to\U(H_s)$,
$$
i^{(n,k)}_X(e_{xy})
 ={1\over|X|}W(Q_nU^kx)W(Q_nU^ky)^*
 ={1\over|X|}e^{-{i\over2}\Ip(Q_nU^kx,Q_nU^ky)}W(Q_nU^k(x-y)).
$$
On the other hand, we have
$$
(\alpha^k_U\circ i_X)(e_{xy})={1\over|X|}W(U^kx)W(U^ky)^*
 ={1\over|X|}e^{-{i\over2}\Ip(U^kx,U^ky)}W(U^k(x-y)).
$$
We may conclude that there exists an upper bound for
$||\alpha^k_U\circ i_X-i^{(n,k)}_X||_{\omega_A}$ depending only on
$\eps$, $||A||$, $|X|$ and $||X||=\max\{||x||\,|\,x\in X\}$. Set
$$
\gamma_{n,k}=i^{(n,k)}_X\otimes(\alpha^k_U\circ i_Z).
$$
Then $||\alpha^k_U\circ\gamma-\gamma_{n,k}||_{\omega_A}$ is bounded by a
value depending only on
$\eps$, $||A||$, $||X||$, $|X|$ and $|Z|$. By Proposition~IV.3 in \cite{CNT},
\begin{equation} \label{e3.2}
|H_{\omega_A}(\gamma,\alpha_U\circ\gamma,\ldots,\alpha^{n-1}_U\circ\gamma)
 -H_{\omega_A}(\gamma_{n,0},\gamma_{n,1},\ldots,\gamma_{n,n-1})|<n\delta,
\end{equation}
where $\delta=\delta(\eps,||A||,||X||,|X|,|Z|)\tends{\eps\to0}0$. Since
$\gamma_{n,k}$'s are channels in $\U(Q_nH_s\oplus H_{m+n-1})$, we have
\begin{equation} \label{e3.3}
H_{\omega_A}(\gamma_{n,0},\gamma_{n,1},\ldots,\gamma_{n,n-1})
 \le S(\omega_A|_{\U(Q_nH_s\oplus H_{m+n-1})}).
\end{equation}
By Lemma~\ref{2.3},
\begin{equation} \label{e3.4}
S(\omega_A|_{\U(Q_nH_s\oplus H_{m+n-1})})
 \le S(\omega_A|_{\U(Q_nH_s)})+S(\omega_A|_{\U(H_{m+n-1})})
\end{equation}
and
\begin{equation} \label{e3.5}
S(\omega_A|_{\U(Q_nH_s)})\le(\Eccr{||A||})\dim Q_nH_s\le\eps n(\Eccr{||A||}).
\end{equation}
>From (\ref{e3.2})-(\ref{e3.5}) we conclude that
$$
h_{\omega_A}(\gamma;\alpha_U)\le\delta+\eps(\Eccr{||A||})
 +\lim_{n\to\infty}{1\over n}S(\omega_A|_{\U(H_{n-1})}).
$$
Because of the arbitrariness of $\eps$, the proof of (\ref{e3.1}) is complete.

Applying Lemma~\ref{2.3}, we obtain

\CAR $\displaystyle h_{\omega_A}(\alpha_U)\le S(\omega_A|_{\A(H_0)})
 =\sum^{m_0}_{m=1}(\Ecar{\lambda_m})$,

\CCR $\displaystyle h_{\omega_A}(\alpha_U)\le S(\omega_A|_{\U(H_0)})
 =\sum^{m_0}_{m=1}(\Eccr{\lambda_m})$,

\noindent
where $\displaystyle \lambda_m=\int_\T\lambda_m(z)d\lambda(z)$. Applying
these inequalities to the operator $U^n$ and using the equality
$h_{\omega_A}(\alpha_U)={1\over n}h_{\omega_A}(\alpha_{U^n})$, we may
conclude that

\CAR $\displaystyle
h_{\omega_A}(\alpha_U)\le{1\over n}\sum^{m_0}_{m=1}\sum^n_{k=1}
  (\Ecar{\lambda_{mnk}}),$

\CCR $\displaystyle
h_{\omega_A}(\alpha_U)\le{1\over n}\sum^{m_0}_{m=1}\sum^n_{k=1}
  (\Eccr{\lambda_{mnk}}),$

\noindent
where $\displaystyle \lambda_{mnk}=n\int^{k\over n}_{k-1\over n}
  \lambda_m\left(e^{2\pi i t}\right)dt$.

It remains to make use of the following lemma.

\begin{lemma} \label{3.1}
 Let $g$ be a bounded measurable function, $f$ a continuous function. Then
 $$
 \lim_{n\to\infty}{1\over n}\sum^n_{k=1}f
   \left(n\int^{k\over n}_{k-1\over n}g(t)dt\right)=\int^1_0f(g(t))dt.
 $$
\end{lemma}

\noindent{\it Proof.} Define a linear operator $F_n$ on $L^1(0,1)$,
$$
(F_nh)(t)=n\int^{k\over n}_{k-1\over n}h(t)dt\ \ \hbox{on}\ \
 \left[{k-1\over n},{k\over n}\right].
$$
Then $F_n\to{\rm id}$ pointwise-norm. Indeed, since $||F_n||=1$, it
suffices to prove the assertion for continuous functions, for which it
is obvious. Thus, $F_ng\to g$ in mean, hence in measure.
By virtue of the uniform continuity of $f$, we conclude that
$f\circ F_ng\to f\circ g$ in measure, whence
$\displaystyle\int^1_0f\circ F_ng\,dt\to\int^1_0f\circ g\,dt$.
\endverif

\bigskip\bigskip
\section{Lower bound for the entropy: basic estimate} \label{4}
The aim of this section is to prove the following estimate.

\begin{proposition} \label{4.1}
 For given $\eps>0$ and $C>0$ ($C<1$ for CAR), there exists $\delta>0$ such
 that if $\Spec A\subset(\lambda_0-\delta,\lambda_0+\delta)$ for some
 $\lambda_0\in(0,C)$ and the spectrum of $U^n$ has Lebesgue component for
 some $n\in\N$, then

 \CAR $\displaystyle
 h_{\omega_A}(\alpha_U|_{\A(H)_e})\ge{1\over n}(\Ecar{\lambda_0}-\eps);$

 \CCR $\displaystyle
 h_{\omega_A}(\alpha_U)\ge{1\over n}(\Eccr{\lambda_0}-\eps).$
\end{proposition}

First, we will prove that if $f\in H$ is close to be an eigenvector for $A$,
then, for any $a\in\CC(\C f)$, $\omega_A(ax)$ is close to
$\omega_A(a)\omega_A(x)$ uniformly on $x\in\CC(\C f)'\cap\CC(H)$.

\begin{lemma} \label{4.2}
 Let $\{e_{ij}\}_{i,j}$ be a system of matrix units in a W$^*$-algebra $M$,
 $e=\sum_ke_{kk}$, $\omega$ a normal faithful state on $M$. Then, for any
 $x\in M$ commuting with the matrix units, we have
 $$
 |\omega(e_{kk}x)-\lambda_k\omega(x)|\le2(\lambda^\2_k||1-e||_\omega+
  \sum_j||\lambda^\2_j\msigma(e_{kj})-\lambda^\2_ke_{kj}||_\omega)
  ||x||^{\#}_\omega,
 $$
 where $||x||^{\#}_\omega=(\omega(x^*x)+\omega(xx^*))^\2$ and
 $\lambda_k=\omega(e_{kk})$.
\end{lemma}

\noindent
{\it Proof.} Let $\xi=\xi_\omega$ and $J=J_\omega$ be the cyclic vector and
the modular involution corresponding to $\omega$. We have

\noindent$\lambda_j\omega(e_{kk}x)=$
$$
=\lambda^\2_j((\lambda^\2_jJe_{jk}-\lambda^\2_ke_{kj})\xi,
  Je_{jk}x\xi)+\lambda^\2_k(e_{kj}Jx^*\xi,(\lambda^\2_jJe_{jk}-
  \lambda^\2_ke_{kj})\xi)+\lambda_k(x\xi,Je_{jj}\xi),
$$
whence
$$
|\lambda_j\omega(e_{kk}x)-\lambda_k(x\xi,Je_{jj}\xi)|\le
 (\lambda^\2_j||x||_\omega+\lambda^\2_k||x^*||_\omega)
  ||\lambda^\2_j\msigma(e_{kj})-\lambda^\2_ke_{kj}||_\omega
$$
\begin{equation} \label{e4.1}
\le2||x||^\#_\omega||\lambda^\2_j\msigma(e_{kj})-\lambda^\2_ke_{kj}||_\omega.
\end{equation}
Further,
\begin{equation} \label{e4.2}
|\omega(e_{kk}x)-\sum_j\lambda_j\omega(e_{kk}x)|=\omega(1-e)|\omega(e_{kk}x)|
 \le||1-e||_\omega\lambda^\2_k||x||_\omega,
\end{equation}
and
\begin{equation} \label{e4.3}
|\sum_j\lambda_k(x\xi,Je_{jj}\xi)-\lambda_k\omega(x)|
 =\lambda_k|(x\xi,J(1-e)\xi)|\le\lambda^\2_k||1-e||_\omega||x||_\omega.
\end{equation}

Summing up (\ref{e4.1})-(\ref{e4.3}), we obtain the desired estimate.
\endverif

Recall that in Section \ref{2} we introduced a system of matrix units
$\{e_{ij}(f)\}_{i,j}$ in $\pi_{\omega_A}(\CC(H))''$ ($f\in H$, $||f||=1$).
In the sequel we will identify $\CC(H)$ with its image in $B(H_{\omega_A})$.

\begin{lemma} \label{4.3}
 \mbox{\ }\newline
 \CAR For given $\eps>0$, there exists $\delta>0$ such that if
 $\Spec A\subset(0,1)$ and
 $$
 \left\|\left({A\over1-A}\right)^\2f
  -\left({\lambda\over1-\lambda}\right)^\2f\right\|<\delta\ \
 \hbox{for some}\ f,\ ||f||=1,\ \hbox{where}\ \lambda=(Af,f),
 $$
 then $\displaystyle
 ||\lambda^\2_j\msigma(e_{kj}(f))-\lambda^\2_ke_{kj}(f)||_{\omega_A}
  \le\eps(\lambda_j\lambda_k)^{1/4},\ k,j=1,2$, where $\lambda_1=1-\lambda$,
 $\lambda_2=\lambda$.

 \noindent
 \CCR For given $\eps>0$, $C>0$ and $k,j\in\Z_+$, there exists $\delta>0$ such
 that if $\Spec A\subset(0,C)$ and
 $$
 \left\|\left({A\over1+A}\right)^\2f
  -\left({\lambda\over1+\lambda}\right)^\2f\right\|<\delta\ \
 \hbox{for some}\ f,\ ||f||=1,\ \hbox{where}\ \lambda=(Af,f),
 $$
 then $\displaystyle
 ||\lambda^\2_j\msigma(e_{kj}(f))-\lambda^\2_ke_{kj}(f)||_{\omega_A}
  \le\eps(\lambda_j\lambda_k)^{1/4}$, where
 $\displaystyle \lambda_m={\lambda^m\over(1+\lambda)^{m+1}}$.
\end{lemma}

\noindent{\it Proof.} We have
$$
||\lambda^\2_j\msigma(e_{kj})-\lambda^\2_ke_{kj}||^2_{\omega_A}
 =2(\lambda_j\lambda_k)^\2((\lambda_j\lambda_k)^\2
   -\omega_A(e_{jk}\msigma(e_{kj}))).
$$
So we must prove that $\omega_A(e_{jk}\msigma(e_{kj}))$ is close to
$(\lambda_j\lambda_k)^\2$ when $\delta$ is sufficiently small.

\CAR Set $\displaystyle B={A\over1-A}$ and
$\displaystyle \beta={\lambda\over1-\lambda}$. We have
$$
\omega_A(e_{12}\msigma(e_{21}))=\omega_A(e_{21}\msigma(e_{12})),
$$
$$
\lambda_1-\omega_A(e_{11}\msigma(e_{11}))=\omega_A(e_{11}\msigma(e_{22}))
 =\lambda_2-\omega_A(e_{22}\msigma(e_{22})).
$$
By virtue of (\ref{e2.1}) and (\ref{e2.2}),
$\msigma(e_{21})=\msigma(a^*(f))=a^*(B^\2f)$, so
$$
\|\msigma(e_{21})-\beta^\2e_{21}\|=\|B^\2f-\beta^\2f\|<\delta,
$$
whence
$$
|\omega_A(e_{12}\msigma(e_{21}))-\lambda^\2(1-\lambda)^\2|
 =|\omega_A(e_{12}(\msigma(e_{21})-\beta^\2e_{21}))|<\delta
$$
and
$$
|\omega_A(e_{11}\msigma(e_{22}))|=|\omega_A(e_{11}(\msigma(e_{21})
  -\beta^\2e_{21})\msigma(e_{12}))|<\delta.
$$

\smallskip
\CCR Set $\displaystyle B={A\over1+A}$ and
$\displaystyle \beta={\lambda\over1+\lambda}$. First consider the case $k=j$.
We have to prove that
$$
\lambda_k-\omega_A(e_{kk}\msigma(e_{kk}))=\omega_A(e_{kk}\msigma(1-e_{kk}))
=\sum_{m\ne k}\omega_A(e_{kk}\msigma(e_{mm}))
$$
is small if $\delta$ is small enough. Since
$\displaystyle
\left\|\sum^\infty_{m=m_0}e_{mm}\right\|_{\omega_A}=\beta^{m_0\over2}
  \le\left(C\over1+C\right)^{m_0\over2}\tends{m_0\to\infty}0$,
it suffices to prove that $\omega_A(e_{kk}\msigma(e_{mm}))$ can be made
arbitrary small for any fixed $m\ne k$. Since
$\omega_A(e_{kk}\msigma(e_{mm}))=\omega_A(e_{mm}\msigma(e_{kk}))$, we may
suppose that $m>k$, i.~e., $m=k+n$ for some $n\in\N$. We have
(see (\ref{e2.5}))
$$
e_{k+n,k+n}=c_{kn}a^*(f)^{k+1}e_{n-1,k+n},\ \ \hbox{where}\ \
 c_{kn}=\left({(n-1)!\over(k+n)!}\right)^\2.
$$
Using (\ref{e2.4}), we obtain
$$
\Delta^\2e_{k+n,k+n}\xi=c_{kn}Je_{k+n,n-1}Ja^*(B^\2f)^{k+1}\xi.
$$
Since $||(a^*(f_1)^{k+1}-a^*(f_2)^{k+1})\xi||$ is bounded by a value which
depends only on $k$, $||A||$, $||f_i||$ and $||f_1-f_2||$ (this is most easily
seen from the explicit description of the GNS-representation in terms of
the Fock representation, see Example 5.2.18 in \cite{BR2}), we conclude that
$\msigma(e_{k+n,k+n})\xi$ is close to
$$
\beta^{k+1\over2}c_{kn}Je_{k+n,n-1}Ja^*(f)^{k+1}\xi
$$
when $\delta$ is sufficiently small. But then $e_{kk}\msigma(e_{k+n,k+n})\xi$
is close to
$$
\beta^{k+1\over2}c_{kn}Je_{k+n,n-1}Je_{kk}a^*(f)^{k+1}\xi=0.
$$

It remains to consider the case $j\ne k$. As above, we may suppose that
$j>k$, $j=k+n$ for some $n\in\N$. We have
$$
e_{k+n,k}=d_{kn}a^*(f)^ne_{kk},\ \ \hbox{where}\ \
 d_{kn}=\left({k!\over(k+n)!}\right)^\2.
$$
As above, we conclude that $\msigma(e_{k+n,k})\xi$ is close to
$$
\beta^{n\over2}d_{kn}Je_{kk}Ja^*(f)^n\xi
$$
for sufficiently small $\delta$, so $\omega_A(e_{k,k+n}\msigma(e_{k+n,k}))$
is close to
$$
\beta^{n\over2}d_{kn}(e_{k,k+n}a^*(f)^n\xi,Je_{kk}\xi)
 =\beta^{n\over2}(e_{kk}\xi,Je_{kk}\xi)
 =\beta^{n\over2}\omega_A(e_{kk}\msigma(e_{kk})).
$$
As we have proved, $\omega_A(e_{kk}\msigma(e_{kk}))$ can be made close to
$\displaystyle\lambda_k={\lambda^k\over(1+\lambda)^{k+1}}$, but then
$\displaystyle \beta^{n\over2}\omega_A(e_{kk}\msigma(e_{kk}))$
is close to
$\displaystyle \left({\lambda\over1+\lambda}\right)^{n\over2}\cdot
  {\lambda^k\over(1+\lambda)^{k+1}}=(\lambda_k\lambda_{k+n})^\2$.
\endverif

\begin{lemma} \label{4.4}
For given $N\in\N$ and $\eps>0$, there exists $\delta=\delta(\eps,N)>0$
such that if $A$ is an abelian W$^*$-algebra, $\omega$ a normal faithful
state on $A$, $B\subset A$ a W$^*$-subalgebra, and $\{x_i\}^N_{i=1}$ a
family of projections in $A$ such that $\sum_ix_i=1$ and
$$
|\omega(x_iy)-\omega(x_i)\omega(y)|\le\delta||y|| \ \ \forall y\in B,\
 i=1,\ldots,N,
$$
then
$$
\sum_{i,j}\eta(\omega(x_iy_j))\ge\sum_i\eta(\omega(x_i))+
   \sum_j\eta(\omega(y_j))-\eps
$$
for any finite family of projections $\{y_j\}_j$ in $B$ with $\sum_jy_j=1$.
\end{lemma}

\noindent{\it Proof.} Cf. \cite[Lemma 3.2]{GN1}.
\endverif

The proof of Theorem 3.1 in \cite{GN1} shows that Lemma~\ref{4.4} is also
valid for non-abelian $A$ (with $x_i\in B'\cap A$) and without the
requirement that $x_i$'s and $y_j$'s are projections, but we will not use
this fact.

\noindent{\it Proof of Proposition \ref{4.1}.}
Consider the case of CCR-algebra.

There exists $\delta_1>0$ such that
$$
|\eta(\lambda)-\eta(1+\lambda)-\eta(\lambda_0)+\eta(1+\lambda_0)|
 <{\eps\over6}\ \ \forall\lambda_0\in(0,C) \
  \forall\lambda\ge0:|\lambda-\lambda_0|<\delta_1.
$$
We can find $N\in\N$ such that
$$
\sum^\infty_{k=N}\eta\left({\lambda^k\over(1+\lambda)^{k+1}}\right)
 <{\eps\over6}\ \ \forall\lambda\in(0,C+\delta_1).
$$
Then, since
$\sum^\infty_{k=0}\eta(\lambda^k(1+\lambda)^{-k-1})
 =\Eccr{\lambda}$, we have
\begin{equation} \label{e4.4}
\sum^{N-1}_{k=0}\eta\left({\lambda^k\over(1+\lambda)^{k+1}}\right)
 >\Eccr{\lambda_0}-{\eps\over3}\ \ \forall\lambda_0\in(0,C) \
 \forall\lambda\ge0:|\lambda-\lambda_0|<\delta_1.
\end{equation}

By assumptions of Proposition, there exists $f\in H$ such that
$\{U^{kn}f\}_{k\in\Z}$ is an orthonormal system in $H$. Set $p_k=e_{kk}(f)$,
$k=0,\ldots,N-1$, and $p_N=1-\sum^{N-1}_{k=0}p_k$. Let $\P$ be the algebra
generated by $p_k$, $k=0,\ldots,N$. Then
$$
h_{\omega_A}(\alpha_U)
 \ge\lim_{k\to\infty}{1\over kn}
       H_{\omega_A}(\P,\alpha_U(\P),\ldots,\alpha^{kn-1}_U(\P))
 \ge\lim_{k\to\infty}{1\over kn}H_{\omega_A}(\P,\alpha^n_U(\P),\ldots,
         \alpha^{k(n-1)}_U(\P))
$$
\hspace*{1.45cm}
$\displaystyle
 \ge\lim_{k\to\infty}{1\over kn}\sum^N_{i_0,\ldots,i_{k-1}=0}
       \eta(\omega_A(p_{i_0}\alpha^n_U(p_{i_1})
          \ldots\alpha^{(k-1)n}_U(p_{i_{k-1}})))+
$
\begin{equation} \label{e4.5}
 +{1\over n}\sum^N_{j=0}S(\omega_A|_\P,\omega_A(\cdot\,\msigma(p_j))|_\P).
\end{equation}
We want to prove that if $\Spec A\subset(\lambda_0-\delta,\lambda_0+\delta)$
with sufficiently small $\delta$, then the first term in (\ref{e4.5}) is
close
to ${1\over n}(\Eccr{\lambda_0})$ to within ${\eps\over n}$, while the second
term is close to zero.

Start with the second term. We have
\begin{eqnarray*}
\sum^N_{j=0}S(\omega_A|_\P,\omega_A(\cdot\,\msigma(p_j))|_\P)
 &=&\sum^N_{j=0}\sum^N_{k=0}\omega_A(p_k\msigma(p_j))
      (\log\omega_A(p_k\msigma(p_j))-\log\omega_A(p_k))\\
 &=&\sum^N_{k=0}\left(\eta(\omega_A(p_k))
      -\sum^N_{j=0}\eta(\omega_A(p_k\msigma(p_j)))\right).
\end{eqnarray*}
By Lemma~\ref{4.3}, $\omega_A(p_k\msigma(p_j))$ can be made arbitrary close
to $\delta_{kj}\omega_A(p_k)$ (more precisely, we can state that this is true
for $j\le N-1$, but since
$\omega_A(p_k\msigma(p_N))=\omega_A(p_N\msigma(p_k))$ and
$$
\omega_A(p_N)-\omega_A(p_N\msigma(p_N))
 =\sum^{N-1}_{k=0}\omega_A(p_N\msigma(p_k)),
$$
this holds for all $k,j\le N$). Hence, there exists $\delta_2\in(0,\delta_1)$
such that if $\Spec A\subset(\lambda_0-\delta_2,\lambda_0+\delta_2)$,
$\lambda_0\in(0,C)$, then
\begin{equation} \label{e4.6}
\sum^N_{j=0}S(\omega_A|_\P,\omega_A(\cdot\,\msigma(p_j))|_\P)>-{\eps\over3}.
\end{equation}

Turning to the first term in (\ref{e4.5}), set
\begin{equation} \label{e4.7}
\eps_1=\delta\left({\eps\over3},N+1\right),
\end{equation}
where $\delta(\cdot\,,\,\cdot)$ is from Lemma~\ref{4.4}. Find
$N_1\in\N$ such that
$$
\left(C+\delta_2\over1+C+\delta_2\right)^{N_1\over2}<{\eps_1\over8N}.
$$
Then
$$
\left\|1-\sum^{N_1-1}_{k=0}e_{kk}(f)\right\|_{\omega_A}<{\eps_1\over8N}\ \
 \hbox{if}\ \ A\le C+\delta_2,
$$
hence, by Lemmas 4.2 and 4.3 applied to $\{e_{kj}(f)\}^{N_1-1}_{k,j=0}$, there
exists $\delta_3\in(0,\delta_2)$ such that
if $\Spec A\subset(\lambda_0-\delta_3,\lambda_0+\delta_3)$,
$\lambda_0\in(0,C)$, then
\begin{equation} \label{e4.8}
|\omega_A(p_kx)-\omega_A(p_k)\omega_A(x)|\le{\eps_1\over2N}||x||^\#_{\omega_A}
 \le\eps_1||x||\ \ \forall x\in\U(f^\perp)'',\ k=0,\ldots,N-1.
\end{equation}
We have also
\begin{equation} \label{e4.9}
|\omega_A(p_Nx)-\omega_A(p_N)\omega_A(x)|\le\sum^{N-1}_{k=0}
|\omega_A(p_kx)-\omega_A(p_k)\omega_A(x)|\le\eps_1||x||\ \
 \forall x\in\U(f^\perp)''.
\end{equation}
>From (\ref{e4.7})-(\ref{e4.9}) and Lemma~\ref{4.4} we infer that
if $\Spec A\subset(\lambda_0-\delta_3,\lambda_0+\delta_3)$,
$\lambda_0\in(0,C)$, then, for any $k\in\N$,

$\displaystyle
\sum^N_{i_0,\ldots,i_{k-1}=0}\eta(\omega_A(p_{i_0}\alpha^n_U(p_{i_1})
  \ldots\alpha^{(k-1)n}_U(p_{i_{k-1}})))\ge$
\begin{eqnarray}
 &\ge&\sum^N_{i_0=0}\eta(\omega_A(p_{i_0}))+
       \sum^N_{i_1,\ldots,i_{k-1}=0}\eta(\omega_A(p_{i_1}\alpha^n_U(p_{i_2})
         \ldots\alpha^{(k-2)n}_U(p_{i_{k-1}})))-{\eps\over3} \nonumber\\
 &\ge&\ldots\ge k\sum^N_{j=0}\eta(\omega_A(p_j))-(k-1){\eps\over3}
       >k\sum^{N-1}_{j=0}\eta\left({\lambda^j\over(1+\lambda)^{j+1}}\right)
         -(k-1){\eps\over3}, \label{e4.10}
\end{eqnarray}
where $\lambda=(Af,f)\in(\lambda_0-\delta_3,\lambda_0+\delta_3)$. It follows
from (\ref{e4.4}), (\ref{e4.6}) and (\ref{e4.10}) that we may take
$\delta=\delta_3$.

The proof for CAR is similar, and we omit the details.
\endverif

\bigskip\bigskip
\section{Lower bound for the entropy: end of the proof} \label{5}
In this section we will complete the proof of the lower bound for the entropy.

By virtue of the existence of an $\omega_A$-preserving conditional
expectation $\CC(H)\to\CC(H_a)$, we have
$h_{\omega_A}(\alpha_U)\ge h_{\omega_A}(\alpha_{U_a})$. So we may suppose
that $U$ has absolutely continuous spectrum.

First, we will extend Proposition~\ref{4.1} to arbitrary unitaries. The main
step here is the following observation.

\begin{lemma} \label{5.1}
 Let $U_n$ be a unitary operator on $H_n$, $n\in\N$, and
 $\{z_n\}^\infty_{n=1}\subset\T$. Consider two unitary operators $U'$ and
 $U''$ on $H=\oplus^\infty_{n=1}H_n$,
 $$
 U'=\mathop{\oplus}^\infty_{n=1}U_n,\ \
  U''=\mathop{\oplus}^\infty_{n=1}z_nU_n.
 $$
 Then $h_{\omega_A}(\alpha_{U'})=h_{\omega_A}(\alpha_{U''})$ for any
 $\alpha_{U'}$- and $\alpha_{U''}$-invariant quasi-free state $\omega_A$ on
 $\CC(H)$. For CAR, the same holds for the restrictions of the automorphisms
 to the even part $\A(H)_e$ of the algebra.
\end{lemma}

\noindent{\it Proof.} For CAR, this was proved in \cite[Lemma~2.4]{GN2}.
For CCR, the result is valid by similar reasons.

Consider the unitary operator $V=\oplus^\infty_{n=1}z_n$. We state that there
exists a set $\{\CC_i\}_i$ of finite-dimensional C$^*$-subalgebras of
$\CC(H)''(\subset B(H_{\omega_A}))$ such that $\alpha_V(\CC_i)=\CC_i$ and
$$
h_{\omega_A}(\alpha)=\sup_ih_{\omega_A}(\CC_i;\alpha)
$$
for any $\omega_A$-preserving automorphism $\alpha$. Suppose the statement
is proved. Then, since $\alpha_{U''}=\alpha_V\alpha_{U'}=\alpha_{U'}\alpha_V$,
we have $\alpha^k_{U''}(\CC_i)=\alpha^k_{U'}(\CC_i)$ $\forall k\in\Z$,
and hence
$h_{\omega_A}(\CC_i;\alpha_{U''})=h_{\omega_A}(\CC_i;\alpha_{U'})$
$\forall i$, whence
$h_{\omega_A}(\alpha_{U''})=h_{\omega_A}(\alpha_{U'})$.

For each $n\in\N$, choose an increasing sequence
$\{H_{nk}\}^\infty_{k=1}$ of finite-dimensional subspaces of $H_n$ such that
$\cup_kH_{nk}$ is dense in $H_n$. Set
$$
K_n=H_{1n}\oplus\ldots\oplus H_{nn}.
$$
Then $K_n$ is finite-dimensional, $K_n\subset K_{n+1}$, $\cup K_n$ is dense
in $H$. Since $VK_n=K_n$, for CAR we may take $\CC_n=\A(K_n)$ (respectively,
for the even part we may take $\A(K_n)_e$). For CCR, we can not take
$\U(K_n)$'s, since they are infinite-dimensional. However there exist
finite-dimensional subalgebras of $\U(K_n)''$ that are still invariant
under $\alpha_V$. Namely,
for any finite-dimensional subspace $K$ of $H$ and any $n\in\N$, we define
a finite-dimensional C$^*$-subalgebra $\U_n(K)$ of
$\U(H)''$ as follows. Let $N_K$ be the number
operator corresponding to $K$, i.~e.,
$$
N_K=a^*(f_1)a(f_1)+\ldots+a^*(f_m)a(f_m),
$$
where $f_1,\ldots,f_m$ is an orthonormal basis in $K$. This is a selfadjoint
operator affiliated with $\U(K)''$, its spectrum is $\Z_+$ (see \cite{BR2}).
Let $P_n(K)$ be the spectral projection of $N_K$ corresponding to $[0,n-1]$.
Set
$$
\U_n(K)=P_n(K)\U(K)''P_n(K)+\C(1-P_n(K)).
$$
The algebra $\U_n(K)$ is finite-dimensional, since in the Fock representation
of $\U(K)$ the projection $P_n(K)$ is the projection onto the first $n$
components of the symmetric Fock space over $K$, and any other regular
representation of $\U(K)$ is quasi-equivalent to the Fock representation.
If $VK=K$, then $\alpha_V(N_K)=N_K$ and $\alpha_V(\U(K))=\U(K)$, hence
$\alpha_V(\U_n(K))=\U_n(K)$. Since $\cup_n\U(K_n)''$ is weakly dense in
$\U(H)''$, and $\cup_m\U_m(K_n)$
is weakly dense in $\U(K_n)''$, by Lemma~\ref{2.2} we conclude that any
channel in $\U(H)''$ can be approximated in strong operator topology by a
channel $\gamma$ in $\U_m(K_n)$ for some $m,n\in\N$. But then
$h_{\omega_A}(\gamma;\alpha)\le h_{\omega_A}(\U_m(K_n);\alpha)$.
Thus we may take $\CC_{mn}=\U_m(K_n)$.
\endverif

\begin{lemma} \label{5.2}
 Let $X_1$, $X_2$ be measurable subsets of $\T$, $\lambda(X_1)$,
 $\lambda(X_2)>0$. Then there exist a measurable subset $Y$ of $X_1$,
 $\lambda(Y)>0$, and $z\in\T$ such that $zY\subset X_2$.
\end{lemma}

\noindent{\it Proof.} See Lemma~3.5 in \cite{GN2} where this lemma is proved
for arbitrary locally compact groups. We want only to note that for $\T$
the result is rather obvious in view of the possibility of approximating
measurable sets by finite unions of arcs.
\endverif

Now we can extend Proposition~4.1 to arbitrary unitaries (with absolutely
continuous spectrum). Consider a direct integral decomposition
$$
H=\int^\oplus_\T H_zd\lambda(z),\ \ U=\int^\oplus_\T z\,d\lambda(z),
$$
and set $X=\{z\in\T\,|\,H_z\ne0\}$.

\begin{lemma} \label{5.3}
 For given $\eps>0$ and $C>0$ ($C<1$ for CAR) there exists $\delta>0$ such
 that if $\Spec A\subset(\lambda_0-\delta,\lambda_0+\delta)$ for some
 $\lambda_0\in(0,C)$, then

 \CAR $\displaystyle
 h_{\omega_A}(\alpha_U|_{\A(H)_e})\ge\lambda(X)(\Ecar{\lambda_0}-\eps);$

 \CCR $\displaystyle
 h_{\omega_A}(\alpha_U)\ge\lambda(X)(\Eccr{\lambda_0}-\eps).$
\end{lemma}

\noindent{\it Proof.} Consider the case of CAR-algebra. Choose $\delta>0$
as in the formulation of Proposition~\ref{4.1}. Let
$\{n_k\}^\infty_{k=1}\subset\N$ be a sequence such that
$\lambda(X)=\sum_k{1\over n_k}$. The Zorn lemma and Lemma~\ref{5.2} ensure
the existence of an at most countable set $\{X_{1m}\}_m$ of disjoint measurable
subsets of $X$ and a set $\{z_{1m}\}_m\subset\T$ such that
\begin{equation} \label{e5.2}
\exp\left(2\pi i\left[0,{1\over n_k}\right]\right)
 =\mathop{\bigsqcup}_mz_{km}X_{km}\ \hbox{mod}\,0
\end{equation}
holds for $k=1$. Proceeding by induction, we obtain a countable measurable
partition $\{X_{km}\}_{k,m}$ of $X$ and a countable subset $\{z_{km}\}_{k,m}$
of $\T$ such that (\ref{e5.2}) holds for all $k\in\N$. Let $H_{km}$ be the
spectral subspace for $U$ corresponding to the set $X_{km}$. Set
$H_k=\oplus_mH_{km}$, and define a unitary operator $U_k$ on $H_k$,
$$
U_k=\mathop{\oplus}_mz_{km}U|_{H_{km}}.
$$
By Lemma~\ref{5.1} and Proposition~\ref{4.1}, we have
$$
h_{\omega_A}(\alpha_U|_{\A(H_k)_e})=h_{\omega_A}(\alpha_{U_k}|_{\A(H_k)_e})
 \ge{1\over n_k}(\Ecar{\lambda_0}-\eps).
$$
For any $k_0\in\N$, there exists an $\omega_A$-preserving conditional
expectation $\A(H)\to\otimes^{k_0}_{k=1}\A(H_k)_e$
(see Remark 4.2 in \cite{SV}). By virtue of the superadditivity of the
entropy \cite[Lemma~3.4]{SV}, we conclude that $$
h_{\omega_A}(\alpha_U|_{\A(H)_e})\ge\sum^{k_0}_{k=1}
   h_{\omega_A}(\alpha_U|_{\A(H_k)_e})\ge\left(\sum^{k_0}_{k=1}
   {1\over n_k}\right)(\Ecar{\lambda_0}-\eps).
$$
Letting $k_0\to\infty$, we obtain the estimate we need.

The proof for CCR is similar, and we omit it.
\endverif

\noindent{\it Proof of Corollary \ref{1.2}.} We will consider only the case of
CAR-algebra. Fix $\delta_0\in(0,{1\over2})$ and take
$\eps\in(0,\eta(\delta_0))$.
Let $\delta$ be as in the formulation of Lemma~\ref{5.3} with $C=1-\delta_0$.
For any Borel subset $X$ of $\R$, let $\1_X(A)$ be the spectral projection
of $A$ corresponding to $X$. Then
$$
\1_X(A)=\int^\oplus_\T\1_X(A_z)d\lambda(z).
$$
Define a measurable function $\phi_X$ on $\T$,
$$
\phi_X(z)=\cases{1,\ \1_X(A_z)\ne0,\cr
                 0,\ \hbox{otherwise}.}
$$
By Lemma~\ref{5.3}, we conclude that if $X$ is a Borel subset of
$(\lambda_0-\delta,\lambda_0+\delta)$ for some
$\lambda_0\in(\delta_0,1-\delta_0)$, then
\begin{equation} \label{e5.3}
h_{\omega_A}(\alpha_U|_{\A(\1_X(A)H)_e})\ge(\Ecar{\lambda_0}-\eps)
   \int_\T\phi_X(z)d\lambda(z)\ge\eta(1-\delta_0)\cdot\int_\T\phi_X(z)d\lambda(z),
\end{equation}
where we have used the inequality $\Ecar{\lambda_0}\ge\Ecar{\delta_0}$.

Let $t_0=\delta_0<t_1<\ldots<t_m=1-\delta_0$, $t_k-t_{k-1}<\delta$. Then by
the same reasons as in the proof of Lemma~\ref{5.3}, we obtain from
(\ref{e5.3}) the inequality
$$
h_{\omega_A}(\alpha_U)\ge\eta(1-\delta_0)\cdot\int_\T\sum^m_{k=1}
   \phi_{(t_{k-1},t_k]}(z)d\lambda(z).
$$
Letting $\max(t_k-t_{k-1})\to0$, we conclude that if
$h_{\omega_A}(\alpha_U) <\infty$, then $(\delta_0,1-\delta_0)\cap\Spec A_z$
is finite for almost all $z\in\T$. Since $\delta_0$ is arbitrary, $A_z$
has pure point for almost all $z$ provided the entropy is finite.
\endverif

It remains to consider the case where $A_z$ has pure point spectrum for
almost all $z$. Then
$$
H=\mathop{\oplus}^N_{n=1}L^2(X_n,d\lambda),
$$
where $X_n$ is a measurable subset of $\T$, $N\le\infty$, $U$ and $A$ act
on $L^2(X_n)$ as multiplications by functions $z$ and $\lambda_n(z)$,
respectively. We must prove that

\CAR $\displaystyle h_{\omega_A}(\alpha_U)\ge\sum^N_{n=1}
   \int_{X_n}(\Ecar{\lambda_n(z)})d\lambda(z)$,

\CCR $\displaystyle h_{\omega_A}(\alpha_U)\ge\sum^N_{n=1}
   \int_{X_n}(\Eccr{\lambda_n(z)})d\lambda(z)$.

Again, consider only the case of CAR-algebra. Using the superadditivity
as above, we see that it suffices to estimate
$h_{\omega_A}(\alpha_U|_{\A(H)_e})$ supposing $N=1$. As in the proof of
Corollary \ref{1.2}, fixing $\delta_0>0$, $\eps>0$ and choosing
$t_0=\delta_0<t_1<\ldots<t_m=1-\delta_0$, we obtain
$$
h_{\omega_A}(\alpha_U|_{\A(H)_e})
 \ge\sum^m_{k=1}\int_{\{t_{k-1}<\lambda_1(z)\le t_k\}}(\Ecar{t_k}-\eps)
     d\lambda(z)
$$
if $\max(t_k-t_{k-1})$ is small enough. Letting $\max(t_k-t_{k-1})\to0$,
we obtain
$$
h_{\omega_A}(\alpha_U|_{\A(H)_e})
 \ge\int_{\{\delta_0<\lambda_1(z)\le1-\delta_0\}}(\Ecar{\lambda_1(z)})
      d\lambda(z)-\eps.
$$
In view of the arbitrariness of $\delta_0$ and $\eps$, the proof is
complete.

\bigskip\bigskip
\section*{Appendix A}
The results of the paper allow to construct a simple example of non-conjugate
K-systems with the same finite entropy (see also Section 5 in \cite{GN1}).

\medskip\noindent
{\bf Theorem A.1} {\it
 Let $U$ be a unitary operator on $H$ with absolutely continuous spectrum,
 $A\in B(H)$, $A\ge0$, $\Ker A=0$, $AU=UA$. Suppose
 $$
 \left({A\over1+A}\right)^{it_0}=U \ \ \hbox{for some}\ \
  t_0\in\R\backslash\{0\}.
 $$
 Let $\omega$ and $\tau_\theta$, $\theta\in\R$, be the quasi-free state and
 the Bogoliubov automorphism of the CCR-algebra $\U(H)$ corresponding to
 $A$ and $e^{i\theta}U$, respectively. Set $M=\pi_\omega(\U(H))''$. Then

 (i) $M$ is the hyperfinite III$_1$-factor;

 (ii) $(M,\omega,\tau_\theta)$, $\theta\in[0,2\pi)$, are pairwise
 non-conjugate entropic K-systems with the same entropy.}

\medskip
\noindent{\it Proof.} There exist a larger space $K\supset H$ and a unitary
operator $V$ on $K$ with homogeneous Lebesgue spectrum such that
$U=V|_H$. Let $C$ be a non-singular bounded positive operator on $K$
commuting with $V$ such that $A=C|_H$. Set $\phi=\omega_C$,
$\beta_\theta=\alpha_{e^{i\theta}V}$ and $N=\pi_\phi(\U(K))''$. Since
$\phi$ is separating, we may consider $M$ as a subalgebra of $N$. The
algebras $M$ and $N$ are hyperfinite III$_1$-factors, moreover, the
centralizer $M_\omega$ is trivial (see, for example, \cite{GN1}, p. 227).
There exists a subspace $K_0$ of $K$ such that $K_0\subset VK_0$,
$\cap_nV^nK_0=0$, $\cup_nV^nK_0$ is dense in $K$. Let $N_0$ be the
W$^*$-subalgebra of $N$ generated by $\U(K_0)$. Then
$N_0\subset\beta_\theta(N_0)$,
$\cup_{n\in\N}(\beta^{-n}_\theta(N_0)'\cap\beta^n_\theta(N_0))\supset
   \cup_{n\in\N}\U(V^nK_0\ominus V^{-n}K_0)$ is weakly dense in $N$,
$\cap_n\beta^n_\theta(N_0)=\C1$ since $N$ is a factor. Hence,
$(N,\phi,\beta_\theta)$ is an entropic K-system by \cite[Theorem 3.1]{GN1}.
Since $(M,\omega,\tau_\theta)$ is a subsystem, and there exists a
$\phi$-preserving conditional expectation $N\to M$, it is an entropic
K-system too.

The fact that $h_\omega(\tau_\theta)$ does not depend on $\theta$ follows
either from the formula for the entropy or directly from Lemma~\ref{5.1}.

It remains to prove the non-conjugacy. Let $\theta\mapsto\gamma_\theta$ be
the gauge action. Since $\tau_\theta=\gamma_\theta\tau_0$, it suffices to
prove that $(M,\omega,\tau_0)$ and $(M,\omega,\tau_\theta)$ are
non-conjugate for $\theta\in(0,2\pi)$. Since $\tau_0=\sigma^\omega_{t_0}$,
any $\omega$-preserving automorphism of $M$ commutes with $\tau_0$ and
can not conjugate $\tau_0$ with an automorphism different from $\tau_0$.
\endverif

Note that any K-automorphism is ergodic, and for any ergodic automorphism
there exists at most one invariant normal state. Hence, any automorphism of
$M$ conjugating $\tau_{\theta_1}$ with $\tau_{\theta_2}$ preserves $\omega$.
Thus, the automorphisms $\tau_\theta$, $\theta\in[0,2\pi)$, are pairwise
non-conjugate (but their restrictions to $\U(H)$ are conjugate).

To obtain finite entropy we may take, for example, unitaries with finitely
multiple spectrum. We see also that if the unitary has homogeneous
Lebesgue spectrum,
then the systems constructed above have the algebraic K-property.

\bigskip\bigskip
\section*{Appendix B}
The following result was used in Sections \ref{3} and \ref{5}.

\smallskip\noindent
{\bf Theorem B.1} {\it Let $(Z,\nu)$ be a standard measure space,
 $Z\ni z\mapsto H_z$
 a measurable field of Hilbert spaces, $d(z)=\dim H_z$,
 $A=\int^\oplus_ZA_zd\nu(z)$ a decomposable selfadjoint operator on
 $H=\int^\oplus_ZH_zd\nu(z)$. Suppose that $A_z$ has pure point spectrum
 $\nu$-a.~e. Then there exist measurable vector fields
 $e_1(z),e_2(z),\ldots\ $, such that $\{e_n(z)\}^{d(z)}_{n=1}$ is an
 orthonormal
 basis of $H_z$ consisting from eigenvectors of $A_z$ for almost all $z$,
 and $e_n(z)=0$ for $n>d(z)$ if $d(z)<\aleph_0$.}
\smallskip

\noindent{\it Proof.} First, prove that there exists a measurable vector
field $e$ such that $e(z)$ is an eigenvector of norm one for $A_z$ for
almost all $z$. In proving this we may suppose that $Z$ is a compact metric
space, $\{H_z\}_z$  the constant field defined by a separable Hilbert space
$H_0$, and $z\mapsto A_z\in B(H_0)$ a weakly continuous mapping.
Consider the subset $X$ of $Z\times H_0\times\R$ defined by
$$
X=\{(z,e,\lambda)\,|\,||e||=1,\, A_ze=\lambda e\}.
$$
Since $X$ is closed, there exists a measurable section for the projection
$X\to Z$, and our statement is proved.

Let $\{e_i\}_{i\in I}$ be a maximal family of vectors in $H$ such that
$e_i(z)$ and $e_j(z)$ are mutually orthogonal a.~e. for $i\ne j$, and
$e_i(z)$ is an eigenvector of norm one for $A_z$ for almost all $z$. Since
$H$ is separable, $I$ is at most countable. Hence,
if $P_z$ is the projection onto the space spanned by $e_i(z),i\in I$, then
$z\mapsto P_z$ is a measurable field of projections, whence
$z\mapsto(1-P_z)H_z$ is a measurable field of subspaces. By the maximality,
$\{e_i(z)\}_{i\in I}$ is an orthonormal basis of $H_z$ consisting from
eigenvectors of $A_z$ on a subset of $Z$ of positive measure. Thus, the
conclusion of Theorem holds on a subset of positive measure. Applying the
maximality argument once again, we obtain an at most countable measurable
partition
of $Z$ such that vector fields with the required properties exist over each
element of the partition. Gluing them, we get the conclusion.
\endverif

Note that if it was a priori known that there exist measurable functions
$\lambda_1(z),\lambda_2(z),\ldots\ $, such that the point spectrum of $A_z$
coincides with $\{\lambda_n(z)\}_n$ (counting with multiplicities),
then the conclusion of Theorem would follow directly from Lemma~2 on p.166
in \cite{D}.

\bigskip

\begin{flushleft}
Institute for Low Temperature Physics \& Engineering\\
Lenin Ave 47\\
Kharkov 310164, Ukraine\\
neshveyev@ilt.kharkov.ua\\
\end{flushleft}

\end{document}